\documentclass[12pt]{article}

\usepackage{amssymb}
\newcommand{\R}{\mathbb R}
\newcommand{\C}{\mathbb C}
\newcommand{\Z}{\mathbb Z}

\renewcommand{\Re}{\mbox{\rm Re }}
\renewcommand{\Im}{\mbox{\rm Im }}
\newcommand{\rz}{\mbox{\rm Res }}

\newtheorem{tetel}{Theorem}
\newtheorem{prop}{Proposition}
\input{tcilatex}

\begin{document}

\title{Bifurcation of Periodic Delay Differential Equations at Points
of 1:4 Resonance \thanks{%
Supported by the Hungarian Foundation for Scientific Research, grant T
049516.} }
\author{Gergely R\"ost \thanks{%
Bolyai Institute, Univ. Szeged, Hungary, H-6720 Szeged, Aradi
v\'ertan\'uk tere 1.} }
\date{}
\maketitle

\begin{abstract}
The time-periodic scalar delay differential equation $\dot
x(t)=\gamma f(t,x(t-1))$ is considered, which leads to a resonant
bifurcation of the equilibrium at critical values of the parameter.
Using Floquet theory, spectral projection and center manifold
reduction, we give conditions for the stability properties of the
bifurcating invariant curves and four-periodic orbits. The
coefficients of the third order normal form are derived explicitly.
We show that the 1:4 resonance has no effect on equations of the
form $\dot z(t)=-\gamma r(t)g(x(t-1))$.
\end{abstract}

\bigskip


\bigskip

\bigskip
{\bf AMS 2000} {34K18, 37M20}

\bigskip

Final version of this paper has been published in Functional Differential Equations, 2006, vol. 13, Nr. 3-4, pp 519-536

\section{Introduction}

The generic bifurcation of planar discrete dynamical systems is the
Neimark-Sacker bifurcation, where a complex conjugate pair of multipliers
crosses the unit circle at critical values of the parameter and an invariant
curve bifurcates from the equilibrium. These results can be extended to
higher dimensional systems by center manifold theorems and projection
methods. The case, when the critical multipliers are fourth roots of unity,
is called 1:4 strong resonance. At strong resonances the Neimark-Sacker
bifurcation theorem is not valid any more and we need further studies. The
modern theory of strong resonances is due to Arnold (\cite{arn88},\cite%
{arn99}). We can not expect the appearance of an invariant curve in general,
and intricate behavior of the system is possible (see e.g. \cite{cheng90},%
\cite{gam}). The most complicated case of resonances is the 1:4. In
this paper we use some related results of Iooss (\cite{iooss}), Wan
(\cite{wan}) and Lemaire (\cite{lemaire}). In a previous paper
(\cite{rg05}) the bifurcation of the time-one map of a scalar
periodic delay differential equation was studied. For the general
case, the bifurcation analysis is performed by using Floquet theory,
center manifold reduction and a spectral projection method, but a
wide class of delay differential equations leads to a strong 1:4
resonance. The purpose of this paper is to broaden the results to
this case.

Consider the non-autonomous scalar delay differential equation

\begin{equation}  \label{egy}
\dot x(t)=\gamma f(t,x(t-1)),
\end{equation}
where $\gamma$ is a real parameter, $f: {\R }\times {\R }\rightarrow
{\R}$ is a $C^4$-smooth function satisfying
\[
f(t+1,\xi)=f(t,\xi)
\]
and
\[
f(t,0)=0
\]
for all $t, \xi \in {\R}.$ Such equations arise very naturally in
several applications, i.e. in population dynamics. A nice overview of
related models can be found in \cite{ruan05}. The periodicity is due to the
periodic fluctuation of the environment. Denote the Banach space of
continuous real and complex valued functions on the interval [-1,0] by $C$
and $C_{{\C}}$, respectively, with the norm
\[
||\phi||=\sup_{-1\leq t \leq 0} |\phi(t)|.
\]
Every $\phi \in C$ determines a unique continuous function $x^\phi:
[-1,\infty) \rightarrow {\R}$, which is differentiable on $(0, \infty)$,
satisfies (\ref{egy}) for all $t>0$ and $x^\phi(t)=\phi(t)$ for all $t \in
[-1,0]$. Such a function $x^\phi$ is called the solution of (1) with the
initial value $\phi$. The time-one map $F: C \rightarrow C$ is defined by
the relations
\[
F(\phi)=x^\phi_1, x_t(s)=x(t+s), s \in [-1,0].
\]
The notation $F_\gamma$ emphasizes the dependence of the time-one map on the
parameter. The spectrum $\sigma(U)$ of the monodromy operator $U$ (the
derivative of the time-one map $F$ at $0$) determines the behavior of
solutions close to the equilibrium $0$. The monodromy operator is a linear
continuous map and with the relation $U(\psi)=U(\Re \psi)+i U(\Im \psi)$
considered as an operator $C_{{\C}} \rightarrow C_{{\C}}$ and given
by $U(\psi)=y_1^\psi$, where $y^\psi$ is the solution of the linear
variational equation

\begin{equation}  \label{lin}
\dot y(t)=\gamma f_\xi(t,0)y(t-1),
\end{equation}
where $y^\psi|_{[-1,0]}\equiv \psi$. The operator $U$ is compact,
therefore all the non-zero points of the spectrum are isolated
points and eigenvalues of finite multiplicity with finite
dimensional range of the associated eigenprojection $P_\mu:C_{{\C}} \rightarrow C_{{\C}}$, where $\mu \in \sigma(U), \mu \neq
0$. These eigenvalues are called Floquet multipliers. The spectral
theory and other properties of different types of delay differential
equations were extensively studied in \cite{deleq} and
\cite{hallun}.

In \cite{rg05}, the equation
\begin{equation}  \label{regi}
\dot x(t)=\gamma \big(a(t)x(t)+f(t,x(t-1))\big)
\end{equation}
was studied. Varying $\gamma$, Floquet multipliers cross the unit
circle and bifurcation of an invariant curve occurs, supposing that
the critical Floquet multipliers are not third or fourth roots of
unity. For equation (\ref{egy}), the cri-tical eigenvalues are $i$
and $-i$, that is a strong resonance, and the results of \cite{rg05}
are not valid anymore. In the case of strong resonance, in general
one can not expect the appearance of the invariant curve (see
\cite{arn88} or \cite{arn99}). In \cite{wan}, a condition was given
which guarantees the appearance of the invariant curve for
2-dimensional maps even at points of resonance. Independently, a
similar result was presented in \cite{lemaire}. Roughly speaking, if
there are no bifurcating four-periodic points, the invariant curve
occurs. The stability of the four-periodic points was treated in
\cite{iooss}. To apply this to our infinite-dimensional system, we
use center manifold reduction.

The classical process of computing the dynamical system restricted
to the center manifold using bilinear forms for delay differential
equations (see \cite{hallun} for the theory and \cite{hkw} for
applications) can not be
applied directly to periodic equations. Faria (\cite{far97} and \cite{far98}%
) presented the method of normal forms for periodic functional
differential equations with autonomous linear part. We established a
spectral projection method in \cite{rg05} for periodic scalar
equations. The spectral projection is represented by a Riesz-Dunford
integral. The resolvent of the monodromy operator of a periodic
delay differential equation and corresponding spectral projections
were calculated in the paper of Frasson and Verduyn Lunel
(\cite[Section 6.2.]{fravl}) in a more general setting. Certain
computations done in \cite{rg05} for equation (\ref{regi}), can be
used for equation (\ref{egy}), simply taking $a(t)\equiv 0$. We
remark that these arguments work only if the period and the delay
are the same. If the delay is not a multiple of the period, then we
can not compute the Floquet multipliers by the characteristic
equation. Some information can be obtained on the Floquet
multipliers in a similar problem in \cite{walsku}, there the period
is three and the delay is one. The most difficult case, when the
delay is incommensurable with the period, there are no results in
this direction.

The paper is organized as follows. In {\it Section 2} we summarize
some previous results, follow by the general theory
(\cite{deleq},\cite{hallun}) or obtained in \cite{rg05}. {\it
Section 3} is devoted to the bifurcation analysis of strong
resonance. We give an explicit condition in terms of $f$ and its
partial derivatives to ensure the bifurcation of an invariant curve
or four-periodic points, and determine the direction of the
appearance and the stability properties. We apply our results to
equations with periodic coefficient in {\it Section 4}, showing that
the resonance does not cause any "anomalies" for this class of
equations. In {\it Section 5} we illustrate the results on the
example of the celebrated Wright equation with periodic coefficient.

\section{Preliminary results}

A non-zero point $\mu$ of the spectrum of the monodromy operator U is called
a Floquet multiplier of equation (\ref{lin}) and any $\lambda$ for which $%
\mu=e^\lambda$ is called a Floquet exponent of equation (\ref{lin}). By the
Floquet theory (\cite[p. 237]{hallun}), $\mu=e^\lambda$ is a Floquet
multiplier of equation (\ref{lin}) if and only if there is a nonzero
solution of equation (\ref{lin}) of the form $y(t)=p(t)e^{\lambda t},$ where
$p(t+1)=p(t)$. Substituting this solution into equation (\ref{lin}), one can
easily deduce that the Floquet exponents are the zeros of the characteristic
function

\begin{equation}  \label{kar}
h(\lambda)= \lambda - \gamma\beta e^{-\lambda},
\end{equation}
where
\[
\beta=\int_{-1}^0f_\xi(t,0)dt.
\]
We assume that $\beta \neq 0$. The eigenfunctions have the form
\[
\chi_\mu(t): [-1,0]\ni t \mapsto e^{\gamma e^{-\lambda}\int_{-1}^t{\
f_\xi(s,0)ds}}\in {\C}.
\]
For any root of the characteristic equation $h(\lambda)=0$, the
corresponding $\chi_\mu(t)$ defines a Floquet solution of equation (\ref{lin}%
), hence the Floquet exponents coincide with the roots of the characteristic
function .

Let
\[
\Delta(z)=z-e^{{\frac{\gamma \beta }{z}}}.
\]
The equation $\Delta(z)=0$ is equivalent to the characteristic
equation. Any complex number $\mu=e^\lambda$ is a root of
$\Delta(z)$ if and only if $ \lambda$ is a Floquet exponent.
Applying Theorem $3.1.$ of \cite[p. 247]{hallun} to equation
(\ref{egy}), one finds that the Floquet multipliers consist of the
roots of $\Delta(z)$ and the algebraic multiplicity of an eigenvalue
$\mu$ equals to the order of $\mu$ as a zero of $\Delta(z)$. When
this number is $1$, we call $\mu$ a simple eigenvalue. According to
the Riesz-Schauder Theorem, if $U: C_{{\C}} \rightarrow C_{{\C}}$ is a compact operator with a simple eigenvalue $\mu$, then
there are two closed subspaces $E_\mu$ and $Q_\mu $ such that $E_\mu
$ is one-dimensional, $E_\mu \oplus Q_\mu =C_{{\C}}$,
furthermore the relations
 $U(E_\mu ) \subset E_\mu $ and $U(Q_\mu ) \subset Q_\mu, \sigma(U|E_\mu)=\{ \mu \}$ and $\sigma(U|Q_\mu)=%
\sigma(U)\backslash \{\mu\}$ hold. The spectral projection $P_\mu$
onto $E_\mu$ along $Q_\mu$ can be represented by the Riesz-Dunford
integral
\[
P_\mu={\frac{1 }{2\pi i }}\int_{\Gamma_\mu} (zI-U)^{-1} dz = \mathop{\rz}%
_{z=\mu} (zI-U)^{-1} ,
\]
where $\Gamma_\mu$ is a small circle around $\mu$ such that $\mu$ is
the only singularity of $(zI-U)^{-1}$ inside $\Gamma_\mu$.

For simplicity, let $b(t)=\gamma f_\xi(t,0)$ and $B(t)=\int_{-1}^t b(s) ds$.
With this notation the linearized equation takes the form
\[
\dot y(t)=b(t)y(t-1),
\]
$\beta={\frac{1 }{\gamma}}\int_{-1}^0b(t)dt$. By the variation-of-constants
formula for ordinary differential equations we find the following
representation of the time-one map

\begin{equation}  \label{F}
F(\phi)(t)=\phi(0)+\int_{-1}^t \gamma f(s,\phi(s)) ds, \quad t \in [-1,0],
\end{equation}
which implies for the monodromy operator

\begin{equation}  \label{U}
U(\phi)(t)=\phi(0)+\int_{-1}^t b(s)\phi(s) ds , \quad t \in [-1,0].
\end{equation}
We need the derivatives of the operator $F$ up to order three, evaluated at $%
0$. Let $V=D^2F(0)$ and $W=D^3F(0)$. $V$ and $W$ are $n$-linear operators
with $n=2$ and $n=3$, respectively. By the representation (\ref{F}), one has

\[
V(\phi _{1},\phi _{2})(t)=\int_{-1}^{t}\gamma f_{\xi \xi }(s,0)\phi
_{1}(s)\phi _{2}(s)ds,\quad t\in \lbrack -1,0],
\]%
and
\[
W(\phi _{1},\phi _{2},\phi _{3})(t)=\int_{-1}^{t}\gamma f_{\xi \xi \xi
}(s,0)\phi _{1}(s)\phi _{2}(s)\phi _{3}(s)ds,\quad t\in \lbrack -1,0].%
\]%
The following statements are special cases of Lemma 4 and Theorem 2 of \cite%
{rg05}, setting $a(t)\equiv 0$.

\begin{prop}
(\cite{rg05}) The resolvent of the monodromy operator can be expressed as

\begin{equation}\label{rez}
(zI-U)^{-1}(\psi )(t)=e^{\int_{-1}^{t}{\frac{b(u)}{z}}du}H\left(
t\right) , \quad t\in \lbrack -1,0],
\end{equation}
where

\bigskip

$H\left( t\right) =\left( {\frac{1}{z}}\psi (0)+e^{\int_{-1}^{0}{\frac{b(u)}{%
z}}du}\int_{-1}^{0}{\frac{1}{z^{2}}}e^{-\int_{-1}^{s}{\frac{b(u)}{z}}%
du}b(s)\psi (s)ds\right) \cdot $

$\cdot \left( \left( z-e^{\int_{-1}^{0}{\frac{b(u)}{z}}du}\right) ^{-1}+{\frac{1}{z%
}}e^{-\int_{-1}^{t}{\frac{b(u)}{z}}du}\psi (t)+\int_{-1}^{t}{\frac{1}{z^{2}}}%
e^{-\int_{-1}^{s}{\frac{b(u)}{z}}du}b(s)\psi (s)ds\right)\!\!.$

\bigskip

\noindent The spectral projection operator, corresponding to a
simple eigenvalue $\mu $, has the representation
\[
P_{\mu }(\psi )=\chi _{\mu }R_{\mu }(\psi ),
\]%
where
\begin{equation}
R_{\mu }(\psi )=\bigg({\frac{1}{\mu +\gamma \beta }}\bigg)\bigg(\psi
(0)+\int_{-1}^{0}{\frac{b(s)\psi (s)}{\chi _{\mu }(s)}}ds\bigg).
\end{equation}
\end{prop}

Notice that $R_{\mu }(\chi _{\mu })=1$. Consider the decomposition
\[
C=T^{c}\oplus T^{su},
\]%
where $T^{c}=\Re E_{\mu }\oplus \Im E_{\mu }$ is the critical
2-dimensional realified center eigenspace corresponding to $\mu $
and spanned by $\{\Re \chi _{\mu },\Im \chi _{\mu }\}$, moreover
$T^{su}=\Re Q_{\mu }\oplus \Im Q_{\mu }$ is the 2-codimensional
realified stable-unstable subspace corresponding to the other part
of $\sigma (U)$. The idea of the projection method is that we
introduce new variables $x,y$ and use them as coordinates on these
subspaces. Suppose we have a map

\[
\tilde x=A_1(x) + g(x,y),
\]
\[
\tilde y=A_2(y)+h(x,y),
\]
where $A_1$ and $A_2$ are linear maps on the corresponding subspaces and
\[
g(0,0)=0,\quad Dg(0,0)=0,
\]
\[
h(0,0)=0,\quad Dh(0,0)=0 .
\]
For $y=M(x)$ we have
\[
\tilde x=A_1(x) + g(x,M(x))
\]
\[
\tilde y=A_2(M(x))+h(x,M(x)).
\]
If $M(x)$ denotes the center manifold then by the invariance $\tilde y=
M(\tilde x)$, and thus
\begin{equation}  \label{wig}
M(A_1(x)+g(x,M(x)))= A_2(M(x))+h(x,M(x)).
\end{equation}
The coefficients of the Taylor-expansion of $M(x)$ can be calculated by this
formula. For details and examples we refer to \cite{kuz} and \cite{wig}. The
computations in the infinite dimensional case can be found in \cite{rg05}.
Represent the Taylor-expansion of $F$ in the form
\[
F(\phi)=U(\phi)+{\frac{1 }{2}} V(\phi,\phi)+ {\frac{1 }{6 }}
W(\phi,\phi,\phi) + O(||\phi||^4).
\]
Let $Z(\phi)= F(\phi) - U(\phi)$ be the nonlinear part of $F$. Now decompose
$\phi \in C$ as
\[
\phi=z\chi_\mu+\bar z \bar \chi_\mu + \psi,
\]
where $z = R_\mu(\phi) \in {\C}$, $z\chi_\mu+\bar z \bar \chi_\mu \in
T^c $ and $\psi \in T^{su}$. The complex variable $z$ is a coordinate on the
2-dimensional real eigenspace $T^c$ and the function $\psi$ is a variable in
$T^{su}$. The subspaces $T^c$ and $T^{su}$ are invariant under $U$. For any
real $\phi$, $\phi \in T^{su}$ if and only if $P_\mu(\phi)=0$. $%
U(\chi_\mu)=\mu\chi_\mu$ implies $U(\bar\chi_\mu)=\bar\mu\bar\chi_\mu$, $%
\overline {R_\mu}=R_{\bar \mu}$.

\begin{prop}
{\rm (\cite{rg05})} The restricted map can be written as

\begin{equation}
\tilde{z}=\mu z+{\frac{1}{2}}\rho _{20}z^{2}+\rho _{11}z\bar{z}+{\frac{1}{2}}%
\rho _{02}{\bar{z}}^{2}+{\frac{1}{2}}\rho _{21}z^{2}\bar{z}+{\frac{1}{6}}%
\rho _{03}\bar{z}^{3}+...,  \label{rest}
\end{equation}%
where

\begin{eqnarray}
\rho _{20} &=&R_{\mu }(V(\chi _{\mu },\chi _{\mu }))  \label{co} \\
\rho _{11} &=&R_{\mu }(V(\chi _{\mu },\bar{\chi}_{\mu }))  \nonumber \\
\rho _{02} &=&R_{\mu }(V(\bar{\chi}_{\mu },\bar{\chi}_{\mu }))  \nonumber \\
\rho _{21} &=&R_{\mu }(W(\chi _{\mu },\chi _{\mu },\bar{\chi}_{\mu
}))+2R_{\mu }(V(\chi _{\mu },(1-U)^{-1}V(\chi _{\mu },\bar{\chi}_{\mu })))+
\nonumber \\
&&+R_{\mu }(V(\bar{\chi}_{\mu },(\mu ^{2}-U)^{-1}V(\chi _{\mu },\chi _{\mu
})))+  \nonumber \\
&&+{\frac{{\frac{1}{\mu }}(1-2\mu )}{1-\mu }}R_{\mu }(V(\chi _{\mu },\chi
_{\mu }))R_{\mu }(V(\chi _{\mu },\bar{\chi}_{\mu }))-  \nonumber \\
&&-{\frac{2}{1-{\frac{1}{\mu }}}}|R_{\mu }(V(\chi _{\mu },\bar{\chi}_{\mu
}))|^{2}-{\frac{\mu }{\mu ^{3}-1}}|R_{\mu }(V(\bar{\chi}_{\mu },\bar{\chi}%
_{\mu }))|^{2},  \nonumber \\
\rho _{03} &=&R_{\mu }(W(\bar{\chi}_{\mu },\bar{\chi}_{\mu },\bar{\chi}_{\mu
}))+3R_{\mu }\Big(V(\bar{\chi}_{\mu },(\mu ^{-2}I-U)^{-1}\cdot  \nonumber \\
&&\cdot \big(V(\bar{\chi}_{\mu },\bar{\chi}_{\mu }-R_{\mu }(V(\bar{\chi}%
_{\mu },\bar{\chi}_{\mu }))\chi _{\mu }-R_{\bar{\mu}}(V(\bar{\chi}_{\mu },%
\bar{\chi}_{\mu }))\bar{\chi}_{\mu }\big)\Big).  \nonumber
\end{eqnarray}
\end{prop}

The coefficients $\rho_{20},\rho_{11},\rho_{02},\rho_{21}$ are computed in
\cite[{\it Section 5}]{rg05}. In the non-resonant case $\rho_{03}$ is not
needed, but can be obtained completely analogously, hence the computation is
omitted here.

\section{Bifurcation of the Time-One Map at Points of Resonance}

Two conditions are formulated in the Neimark-Sacker bifurcation theorem: the
transversality condition, and the non-resonance condition, viz. ${\frac{%
\partial \mu(\gamma) }{{\partial\gamma} }}\mid_{\gamma_{j}} \neq 0$ and $%
\mu_j^3 \neq 1, \mu_j^4\neq 1$, where $\gamma_j$ is a critical parameter
value and $\mu_j$ is a corresponding critical multiplier. The following two
lemmas show that the transversality condition is always fulfilled for
equation (\ref{egy}), while $\mu_j^4=1$. This situation is a 1:4 strong
resonance.

\begin{lemma}
The critical values of (\ref{lin}) are
\[
\gamma_j=\frac{-\frac{\pi}{2}+2j\pi}{\beta}, \quad j \in {\Z},
\]
and the corresponding critical Floquet multipliers are $\mu_j=e^{\lambda_j}
=i$ and $\bar \mu_j = e^{\bar \lambda_j} =-i.$ These Floquet-multipliers are
simple eigenvalues and the critical eigenfunctions are
\[
\chi_{\pm i}(t):[-1,0]\ni t \mapsto e^{\mp iB(t)}\in {\C}.
\]
\end{lemma}

{\sl Proof} One can check easily that $i$ and $-i$ can not be a
double root of $\Delta(z)$, thus if $i$ or $-i$ is a
Floquet-multiplier, then it is always a simple eigenvalue. Suppose
that $\lambda=i \theta$ is a critical
Floquet-exponent, then by the real part of (\ref{kar}) we have $%
\cos(\theta)=0$, hence $\theta=\frac{\pi}{2}+2k\pi$ or $\theta=-\frac{\pi}{2}%
+2k\pi$, where $k \in {\Z}$. Taking into account the imaginary part of (%
\ref{kar}), both options lead to the statements of the lemma by
simple calculations. \rightline{\endproof}

 Introduce the notation $B=B(0)=\gamma \beta$.

\begin{lemma}
\[
{\frac{\partial \mu(\gamma) }{{\partial\gamma} }}\mid_{\gamma_{j}}= \frac{%
\beta}{1+\lambda(\gamma_j)}=\frac{\beta}{{1+B^2}}(1+iB)
\]
\end{lemma}

{\sl Proof} By the characteristic equation and the Implicit Function Theorem
$\mu(\gamma)=e^{\frac{\gamma\beta}{\mu(\gamma)}}$ is defined in a
neighborhood of $\gamma_j$. Differentiating with respect to $\gamma$ gives
\[
\mu^{\prime}(\gamma)=e^{\frac{\gamma\beta}{\mu(\gamma)}}\Big(\frac{\beta
\mu(\gamma)-\beta\gamma \mu^{\prime}(\gamma)}{\mu^2(\gamma)}\Big)%
=\beta-\lambda(\gamma) \mu^{\prime}(\gamma).
\]
This yields $\mu^{\prime}(\gamma)=\frac{\beta}{1+\lambda(\gamma)}$. Setting $%
\gamma=\gamma_j$ one has $\lambda=-i\gamma_j\beta=-iB$ and the lemma is
proved.

\rightline{\endproof}

The next proposition is the Poincar\'{e} normal form map for 1:4 resonance.

\begin{prop}
{\rm (\cite[p. 436]{kuz})} Suppose that we have a map $g=g(\gamma ):{\C}%
\mapsto {\C}$, depending on the parameter $\gamma \in {\R}$, and $g$
has the form
\begin{eqnarray}
g(z) &=&\mu z+\frac{\rho _{20}}{2}z^{2}+\rho _{11}z\bar{z}+\frac{\rho _{02}}{%
2}{\bar{z}}^{2}+\frac{\rho _{30}}{6}z^{3}  \label{map} \\
&&+\frac{\rho _{21}}{2}z^{2}{\bar{z}}+\frac{\rho _{12}}{2}z{\bar{z}}^{2}+%
\frac{\rho _{03}}{6}{\bar{z}}^{3}+{\cal O}(|z|^{4}),  \nonumber
\end{eqnarray}%
where $\mu =\mu (\gamma )$ and $\rho _{kl}=\rho _{kl}(\gamma )$ depends on
the parameter smoothly and $\mu (\gamma _{j})=i$ for some critical value $%
\gamma =\gamma _{j}$. Then by a coordinate transformation depending smoothly
on the parameter, in the critical case the transformed map takes the form

\[
\tilde g(w)=i w+ c_1 w^2 \bar w + c_2 \bar w^3 + {\cal O}(|w|^4),
\]
where

\[
c_{1}=\frac{1+3i}{4}\rho_{20}\rho_{11}+ \frac{1-i}{2}\rho_{11}\bar \rho_{11}+%
\frac{-1-i}{4}\rho_{02}\bar \rho_{02}+\frac{\rho_{21}}{2}
\]
and
\[
c_2=\frac{i-1}{4}\rho_{11} \rho_{02}+ \frac{-i-1}{4} \rho_{02}\bar \rho_{20}+%
\frac{\rho_{03}}{6}.
\]
\end{prop}

Note that similar, but different formulas are presented for $c_2$ in
\cite[ Chapter IV]{iooss} and \cite{wan}. These formulas are
miscalculated and false. One can check directly by a
straightforward, but rather elaborative computation that the formula
of \cite[p. 436]{kuz}, presented in Proposition 3 is the correct
one. However, in the literature the wrong formula of \cite {iooss}
is spreading, see for example the recent papers \cite{xie01} and
\cite{xie05}, where applications of the resonant normal form to
mechanical systems are presented. Since the applied formula is not
correct, the obtained results may not be correct as well.

Define $a_1=\frac{c_1}{i}$ , $a_2=\frac{c_2}{i}$ and $d={\frac{\partial
|\mu(\gamma)| }{{\partial\gamma} }}\mid_{\gamma=\gamma_j}$.

\begin{prop}[Resonant bifurcation theorem,\protect\cite{iooss},\protect\cite%
{wan}]
Suppose that we have a map $g(z):{\C }\mapsto {\C}$ of the form (\ref%
{map}), depending smoothly on the parameter $\gamma$, satisfying $d\neq 0$
and $\mu(\gamma_j)=i$.

If $|\Im(\frac{a_1}{d})|>|\frac{a_2}{d}|$, then a unique invariant curve
bifurcates (and no periodic points of order $4$) from the equilibrium $0$ as
the parameter $\gamma$ passes through $\gamma_j$. The cases $\Re a_1<0$ and $%
\Re a_1>0$ are called supercritical and subcritical bifurcations. In the
supercritical case a stable invariant curve appears for $\gamma > \gamma_j$,
while in the subcritical case an unstable invariant curve disappears when $%
\gamma$ increases through $\gamma_j$.

If $|\Im(\frac{a_1}{d})|<|\frac{a_2}{d}|$, then two families of periodic
points of order $4$ bifurcate (and no invariant curve). Moreover, if $%
|a_1|>|a_2|$, the two families bifurcate on the same side and at least one
of them is unstable. If $|a_1|<|a_2|$, then the two families bifurcate on
opposite sides and both of them are unstable.
\end{prop}

\begin{lemma}
For the restricted map of the time-one map corresponding to equation
(\ref{egy}) we have

\begin{eqnarray*}
a_{1} &=&\frac{3-i}{4}\rho _{20}\rho _{11}-\frac{1+i}{2}|\rho _{11}|^{2}-%
\frac{1-i}{4}|\rho _{02}|^{2}-\frac{i}{2}\rho _{21}= \\
&=&-\frac{i}{2}\Big[R_{i}(W(\chi _{i},\chi _{i},\bar{\chi}%
_{i}))+2R_{i}(V(\chi _{i},(1-U)^{-1}V(\chi _{i},\bar{\chi}_{i})))+ \\
&&+R_{i}(V(\bar{\chi}_{i},(i^{2}-U)^{-1}V(\chi _{i},\chi _{i})))\Big],
\end{eqnarray*}

and

\begin{eqnarray*}
a_{2} &=&\frac{-1+i}{4}\bar{\rho _{20}}\rho _{02}+\frac{1+i}{4}\rho
_{11}\rho _{02}-\frac{i}{6}\rho _{03}= \\
&=&-\frac{i}{6}\Big[R_{i}(W(\bar{\chi}_{i},\bar{\chi}_{i},\bar{\chi}%
_{i}))+3R_{i}\Big(V(\bar{\chi}_{i},(i^{-2}I-U)^{-1}\big(V(\bar{\chi}_{i},%
\bar{\chi}_{i})\big)\Big)\Big].
\end{eqnarray*}
\end{lemma}

{\sl Proof} Apply Proposition 2 and Proposition 3 with $\mu=i$. We
obtain the lemma by a simple calculation.

\rightline{\endproof}

Let us define
\begin{equation}  \label{delta}
\delta= |\Im(a_1)-B\Re(a_1)|- |a_2|\sqrt{1+B^2}.
\end{equation}
Remark that $\delta$ depends on the parameter. Some additional computation
yields

\[
|\Im(\frac{a_1}{d})|>|\frac{a_2}{d}| \Leftrightarrow |\Im({a_1}(1-iB)|>|{a_2}%
(1-iB)|,
\]
that is $\delta >0$.

We apply the center manifold theorem for maps in Banach-spaces to the
time-one map $F$. See \cite{kww} for the existence and \cite{fariawu02} for
the smoothness result. Summarizing all the previous lemmas and propositions
of {\it Section 2} and {\it Section 3}, combining with the center manifold
theorem and the reduction principle (for details see \cite{carr},\cite{kuz}
and \cite{wig}), we obtain our main theorem.

\begin{tetel}
The family of time-one maps $F_\gamma$, corresponding to equation (\ref{egy}%
), has at the critical value $\gamma=\gamma_j$ the fixed point $\phi=0$ with
exactly two simple Floquet-multipliers $\mu_j=i$ and $\bar \mu_j=-i$ on the
unit circle. This is a 1:4 strong resonance. The transversality condition is
fulfilled. There is a neighborhood of $0$ in which a unique invariant curve
(and no 4-periodic points) bifurcates from $0$, providing that $\delta>0.$
The direction of the bifurcation is determined by the sign of $\Re(a_1)$. If
$\delta < 0,$ then two families of 4-periodic points (and no invariant
curve) bifurcate from the equilibrium in a neighborhood of $0$. Furthermore,
if $|a_1|>|a_2|$, the two families bifurcate on the same side and at least
one of them is unstable. If $|a_1|<|a_2|$, then the two families bifurcate
on the opposite side and both of them are unstable.
\end{tetel}

The conditions given in the theorem can be checked for any given equation,
we can compute $\gamma_j,$ $a_1$, $a_2$ and $B$ explicitly by terms of $%
f(t,\xi)$ and its partial derivatives.

\section{Equations with Periodic Coefficient}

In this section we consider the equation
\begin{equation}  \label{pc}
\dot z(t) =- \gamma r(t)g(z(t-1)),
\end{equation}
where $\gamma$ is a real parameter, $r: {\R }\rightarrow {\R}$ is a
continuous function satisfying $r(t+1)=r(t)$ for all $t \in {\R}$, $%
g(\xi)$ is a $C^4$-smooth function satisfying $g(0)=0$. Without loss of
generality we may suppose that
\[
g(\xi)=\xi+\frac{S}{2}\xi^2+\frac{T}{6}\xi^3+{\cal O}(\xi^4),
\]
where $S,T \in {\R}$. With our previous notations we have
\[
f(t,\xi)=- r(t)g(\xi),
\]

\[
f_\xi(t,0)=-r(t),
\]
\[
f_{\xi\xi}(t,0)=-Sr(t),
\]
\[
f_{\xi\xi\xi}(t,0)=-Tr(t),
\]
and
\[
b(t)=-\gamma r(t).
\]
We show that this equation behaves at the bifurcation points as a
nonresonant equation: an invariant curve bifurcates and no $4$-periodic
points from the equilibrium $0$. The following lemma is used many times
during the detailed computations.

\begin{lemma}
Let $B(t)=\int_{-1}^t b(s) ds$. Then
\[
\int_{-1}^t e^{B(s)} b(s) ds=e^{B(t)}-1,
\]
\[
\int_{-1}^t e^{B(s)} b(s)B(s) ds=e^{B(t)}B(t)-e^{B(t)}+1.
\]
\end{lemma}

{\sl Proof} The first identity is obvious, the second can be deduced from
the first by a partial integration.

\rightline{\endproof}

\begin{tetel}
For any family of time-one maps corresponding to equation (\ref{pc}), if $T
\neq\frac{11 S^2}{5}$ then a unique invariant curve bifurcates from the
equilibrium $0$ as the parameter $\gamma$ passes through $\gamma_j$. The
bifurcation is supercritical if $T < S^2\big(\frac{11B+2}{5B}\big)$ and
subcritical if $T > S^2\big(\frac{11B+2}{5B}\big)$.
\end{tetel}

{\sl Proof} Let us fix $\gamma=\gamma_j$ to be a critical parameter
value. Using Lemma 4, we have

\[
V(\chi_i,\chi_i)(t)=\int_{-1}^t Sb(s)e^{-2iB(s)}ds= \frac{S}{-2i}%
(e^{-2iB(t)}-1)=\frac{iS}{2}(e^{-2iB(t)}-1),
\]

\[
V(\chi_i, \bar \chi_i)(t)=SB(t),
\]

\[
V(\bar \chi_i, \bar \chi_i)(t)=\int_{-1}^tS b(s)e^{2iB(s)}ds= \frac{S}{2i}%
(e^{2iB(t)}-1)=\frac{-iS}{2}(e^{2iB(t)}-1),
\]

\[
W(\chi_i,\chi_i,\bar \chi_i)=\int_{-1}^t Tb(s)e^{-iB(s)}ds=iT(e^{-iB(t)}-1),
\]
and

\[
W(\bar \chi_i,\bar \chi_i, \bar \chi_i)=\int_{-1}^t Tb(s)e^{3iB(s)}ds=\frac{%
-iT}{3}(e^{3iB(t)}-1).
\]
Notice that $B=-\frac{\pi}{2}+2j\pi$, hence $e^{iB}=\cos B + i \sin B = -i$,
and $e^{imB}=(-i)^m$ for all $m \in Z$. Taking into account this fact, one
obtains

\begin{eqnarray}
R_{i}(e^{miB(t)}) &=&\left( {\frac{1}{i+B}}\right) \left(
e^{miB}+\int_{-1}^{0}b(s)e^{(m+1)iB(s)}ds\right) =  \label{remi} \\
&=&\left( {\frac{1}{i+B}}\right) \left( e^{miB}+\frac{1}{(m+1)i}%
(e^{(m+1)iB}-1)\right) =  \nonumber \\
&=&\left( {\frac{1}{i+B}}\right) \left( (-i)^{m}-i\frac{(-i)^{m+1}-1}{m+1}%
\right) =  \nonumber \\
&=&\frac{m(-i)^{m}+i}{(i+B)(m+1)}  \nonumber
\end{eqnarray}%
for any $m\neq -1$. Observe that $R_{i}(e^{3iB(t)})=\frac{i}{i+B}=R_{i}(1).$
If $m=-1$, we get the eigenfunction $e^{-iB(t)}$, and as in general,

\begin{eqnarray*}
R_{i}(e^{-iB(t)}) &=&R_{i}(\chi _{i}(t))= \\
&=&\bigg({\frac{1}{i+B}}\bigg)\bigg(e^{-iB}+\int_{-1}^{0}b(s)ds\bigg)={\frac{%
i+B}{i+B}}=1.
\end{eqnarray*}

Now let us evaluate the resolvent by Proposition 1 and Lemma 4:

\begin{eqnarray*}
(1-U)^{-1}V(\chi _{i},\bar{\chi}_{i}) &=&e^{B(t)}\Big((SB+e^{B}%
\int_{-1}^{0}e^{-B(s)}b(s)SB(s)ds)(1-e^{B})^{-1} \\
&&+e^{-B(t)}SB(t)+\int_{-1}^{t}e^{-B(s)}b(s)SB(s)ds\Big) \\
&=&Se^{B(t)}\Big(\big(B+e^{B}(-e^{-B}B-e^{-B}+1)\big)(1-e^{B})^{-1} \\
&&+e^{-B(t)}B(t)-e^{-B(t)}B(t)-e^{-B(t)}+1\Big) \\
&=&Se^{B(t)}(-1-e^{-B(t)}+1)=-S.
\end{eqnarray*}
\newline
Referring to Lemma 3 and $(i^{2}-U)^{-1}=(i^{-2}-U)^{-1}=(-1-U)^{-1}$, we
still need

\label{Last rep}

\begin{equation}
(-1-U)^{-1}(e^{miB(t)})=e^{-B(t)}H_{1}\left( t\right) ,  \label{rez-}
\end{equation}%
where

\bigskip

$H_{1}\left( t\right) =(-1-U)^{-1}(e^{miB(t)})=e^{-B(t)}\Big(\big(%
-e^{miB}+e^{-B}\int_{-1}^{0}b(s)e^{(mi+1)B(s)}ds\big)\cdot $

$\cdot \big(-1-e^{-B}\big)^{-1}-e^{(1+mi)B(t)}+%
\int_{-1}^{t}b(s)e^{(1+mi)B(s)}ds\Big)=$

$=e^{-B(t)}\!\Big(\!\big(\!-e^{miB}\!+e^{-B}\frac{e^{(mi+1)B}-1}{mi+1}\big)\!%
\big(\!-1\!-e^{-B}\big)^{-1}-e^{(1+mi)B(t)}+\frac{e^{(mi+1)B(t)}-1}{mi+1}%
\Big)=$

$=e^{-B(t)}\frac{(-i)^{m}mi-1}{(1+e^{-B})(mi+1)}-e^{miB(t)}\frac{mi}{mi+1}.$

\bigskip

\noindent Particularly,

\bigskip

$(-1-U)^{-1}(e^{2iB(t)}-1)=$

$=e^{-B(t)}\frac{(-i)^{2}2i-1}{(1+e^{-B})(2i+1)}-e^{2iB(t)}\frac{2i}{2i+1}%
-e^{-B(t)}\frac{-1}{1+e^{-B}}=$

$=-e^{2iB(t)}\frac{2i}{2i+1}.$

\bigskip

\noindent Similarly, one finds
\begin{equation}
(-1-U)^{-1}(e^{-2iB(t)}-1)=e^{-2iB(t)}\frac{2i}{-2i+1}.
\end{equation}%
Now we are ready to compute the coefficients of the normal form given in
Lemma 3, namely

\begin{eqnarray}
\quad\quad a_{2} &=&-\frac{i}{6}R_{i}\Big[W(\bar{\chi}_{i},\bar{\chi}_{i},\bar{\chi}%
_{i})+3V(\bar{\chi}_{i},(-1-U)^{-1}\big(V(\bar{\chi}_{i},\bar{\chi}_{i})\big)%
\Big]=  \label{a2} \\
&=&-\frac{i}{6}R_{i}\Big[\frac{-iT}{3}(e^{3iB(t)}-1)-\frac{3iS}{2}V(\bar{\chi%
}_{i},(-1-U)^{-1}\big(e^{2iB(t)}-1\big)\Big]  \nonumber \\
&=&-\frac{T}{18}R_{i}\Big[(e^{3iB(t)}-1)\Big]-\frac{S}{4}R_{i}\Big[V(\bar{%
\chi}_{i},-e^{2iB(t)}\frac{2i}{2i+1}\Big]=0,  \nonumber
\end{eqnarray}
where we used the linearity of $R_{i}$ and the identity
\[
\frac{1}{S}R_{i}(V(\bar{\chi}_{i},e^{2iB(t)}))=\frac{1}{T}R_{i}(W(\bar{\chi}%
_{i},\bar{\chi}_{i}\bar{\chi}_{i})=\frac{1}{3i}R_{i}(e^{3iB(t)}-1)=0.
\]%
We use (\ref{remi}) and (\ref{rez-}) to conclude%
\begin{eqnarray}
\quad\quad a_{1} &=&-\frac{i}{2}R_{i}\Big[W(\chi _{i},\chi
_{i},\bar{\chi}_{i})+2V(\chi
_{i},-S)+  \label{a1} \\
&&+V(\bar{\chi}_{i},(-1-U)^{-1}V(\chi _{i},\chi _{i}))\Big]  \nonumber \\
&=&-\frac{i}{2}R_{i}\Big[Ti(e^{-iB(t)}-1)+2(-i)S^{2}(e^{-iB(t)}-1)\Big]
\nonumber \\
&&-\frac{i}{2}R_{i}\Big[V\big(\bar{\chi}_{i},(-1-U)^{-1}\Big(\frac{iS}{2}%
(e^{-2iB(t)}-1)\big)\Big)\Big]  \nonumber \\
&=&\frac{T-2S^{2}}{2}R_{i}[e^{-iB(t)}-1]+\frac{S}{4}R_{i}\Big[V\big(\bar{\chi%
}_{i},e^{-2iB(t)}\frac{2i}{-2i+1}\Big]  \nonumber \\
&=&\frac{(T-2S^{2})B}{2(i+B)}+\frac{Si}{2-4i}R_{i}[V(\bar{\chi}%
_{i},e^{-2iB(t)})]  \nonumber \\
&=&\frac{(T-2S^{2})B}{2(i+B)}+\frac{Si}{2-4i}R_{i}[\frac{S}{-i}%
(e^{-iB(t)}-1)]  \nonumber \\
&=&\frac{(T-2S^{2})B}{2(i+B)}-\frac{S^{2}B}{(2-4i)(i+B)}=\frac{B}{2(i+B)}%
\big(T-S^{2}\frac{11+2i}{5}\big).  \nonumber
\end{eqnarray}

Applying (\ref{a2}),(\ref{a1}) and $\frac{1}{i+B}=\frac{B-i}{{1+B^{2}}}$ to (%
\ref{delta}), we find
\begin{equation}
2(1+B^{2})\Re (a_{1})=TB^{2}-BS^{2}\frac{11B+2}{5},
\end{equation}

\[
2(1+B^{2})\Im (a_{1})=-TB+BS^{2}\frac{11-2B}{5}.
\]%
The sign of $\Re (a_{1})$ determines the direction of the bifurcation, as
formulated in the theorem, which is the same as the sign of $T-S^{2}\big(%
\frac{11B+2}{5B}\big)$. Substituting the previous two formulas into (\ref%
{delta}), we deduce

\begin{eqnarray*}
\delta &=&\frac{1}{2(1+B^{2})}|B|\cdot |-T-TB^{2}+\frac{S^{2}}{5}%
(11-2B+2B+11B^{2})| \\
&=&\frac{|B|}{2}\cdot |T-S^{2}\frac{11}{5}|.
\end{eqnarray*}%
Since $B\neq 0$, the condition $T\neq \frac{11S^{2}}{5}$ guarantees that $%
\delta >0$ and Theorem 2 is proved.

\rightline{\endproof}

\section{An example}

The classical form of the celebrated Wright-Hutchinson equation (or delayed
logistic equation) is
\[
\dot y(t)=-\alpha y(t-1)(1+y(t)).
\]
The change of variable $z(t)=\ln(1+y(t))$ transforms Wright's equation into
\[
\dot z(t)=-\alpha (e^{z(t-1)}-1).
\]
Since the pioneer works of Wright (\cite{wright55}), a huge amount of papers
concerned with the dynamical properties of this equation and its
generalizations. Here we consider this equation with a periodic coefficient:
\begin{equation}
\dot z(t)=-\alpha r(t)(e^{z(t-1)}-1),  \label{wri}
\end{equation}
where $\alpha >0$ and $r(t)$ is a continuous function satisfying $%
r(t+1)=r(t) $ for all $t \in {\R}$. Without loss of generality we may
suppose $\int_{-1}^0 r(s) ds=1$. We have $g(\xi)=\xi+\frac{1}{2}\xi^2+\frac{1%
}{6}\xi^3+{\cal O}(\xi^4),$ that is $S=1$ and $T=1$. The next theorem is a
direct application of Theorem 2 and Lemma 1.

\begin{tetel}
The family of time-one maps corresponding to equation (\ref{wri}), undergoes
a supercritical bifurcation and a unique invariant curve bifurcates from the
equilibrium $0$ as the parameter $\alpha$ passes through $\frac{\pi}{2}$.
\end{tetel}

Remark that taking $r(t)\equiv 1$ we get back the autonomous case. For the
autonomous Wright equation it is well known that at the value $\alpha=\frac{%
\pi}{2}$ a periodic solution emerges from the equilibrium by a supercritical
Hopf bifurcation. This is consistent with Theorem 3.

\bigskip

{\bf Acknowledgement}.\newline
The author would like to thank the invaluable work of the referees.

\end{document}